# TREE-STRUCTURED REGRESSION AND THE DIFFERENTIATION OF INTEGRALS[1]

By Richard A. Olshen

*Stanford University*

*In Memory of Leo Breiman, 1928–2005*

This paper provides answers to questions regarding the almost sure limiting behavior of rooted, binary tree-structured rules for regression. Examples show that questions raised by Gordon and Olshen in 1984 have negative answers. For these examples of regression functions and sequences of their associated binary tree-structured approximations, for all regression functions except those in a set of the first category, almost sure consistency fails dramatically on events of full probability. One consequence is that almost sure consistency of binary tree-structured rules such as CART requires conditions beyond requiring that (1) the regression function be in $\mathcal{L}^1$, (2) partitions of a Euclidean feature space be into polytopes with sides parallel to coordinate axes, (3) the mesh of the partitions becomes arbitrarily fine almost surely and (4) the empirical learning sample content of each polytope be "large enough." The material in this paper includes the solution to a problem raised by Dudley in discussions. The main results have a corollary regarding the lack of almost sure consistency of certain Bayes-risk consistent rules for classification.

**1. Introduction.** Rooted, binary tree-structured methods have been important modern statistical tools for regression, classification, probability class estimation, clustering and survival analysis; see books by Breiman, Friedman, Olshen and Stone [1], Gersho and Gray [7], Devroye, Györfi and Lugosi [4], Ripley [12], Zhang and Singer [15], Hastie, Tibshirani and Friedman [9] and their references. These books include algorithms, wide ranging

Received November 2004; accepted April 2006.

[1]Supported in part by NSF Grant CCF-00-73050 and in part by NIH/NIBIB Grant 5RO1 EB002784-28.

AMS 2000 subject classifications. 26B05, 28A15, 62G08, 62C12.

Key words and phrases. Binary tree-structured partitions, regression, maximal functions, differentiation of integrals.

applications, and theory. The last has involved an "empirical Lebesgue integral" (an expression first used by Peter Huber), along with connections to the asymptotically minimax approximation of functions (see [6]), and has motivated improvements to the celebrated large deviation result of Vapnik and Chervonenkis; see [10, 11]. To put this paper in context, see [5, 8].

My primary goal is to answer in the negative questions raised by Gordon and Olshen [8] regarding the almost sure limiting behavior of rooted, binary tree-structured rules for regression. There is also solution to a problem posed by Dudley in discussions. The arguments regarding regression can be applied to obtain a certain negative result concerning classification. The remainder of this section is an introduction to terminology and the results in the remainder of the paper. The first part of Section 2 is a summary of relevant results on martingales, on the differentiation of integrals and also on equivariance. It is intended to place Theorems 1.2 and 1.3 and Corollary 1.4 in the somewhat subtle context of previous work. Readers will see that conclusions have two distinct parts, one probabilistic and the other concerning the differentiation of integrals; see (4.1). Lemma 2.1 and Section 3 are expositions of the key parts of the counterexamples; they concern the differentiation of integrals insofar as it is related to rooted, binary tree-structured statistical rules. Section 4 provides details by means of which proofs of the two theorems and the corollary are completed.

As in Breiman, Friedman, Olshen and Stone ([1], Chapter 10), a *finite rooted binary tree* is a finite nonempty set $T$ of positive integers together with (for $t \in T$) two functions, $left(t)$ and $right(t)$, that map $T$ to $T \cup \{0\}$ and which satisfy the following two properties: (1) for each $t \in T$, either $left(t) = right(t) = 0$, or $left(t) > t$ and $right(t) > t$; (2) for each $t \in T$, other than the smallest integer in $T$, there is exactly one $s \in T$ for which either $t = left(s)$ or $t = right(s)$. The minimum element of $T$ is called the *root* of $T$. If $s, t \in T$ and $t = left(s)$ or $t = right(s)$, then $s$ is called the *parent* of $t$. The root of $T$ has no parent, but every other $t \in T$ has a unique parent. A $t \in T$ is called a *terminal node* if it is not a parent, that is, if $left(t) = right(t) = 0$. A finite partition of a set $\Omega$ is called a *finite, rooted, binary tree-structured partition* if there exist a finite, rooted, binary tree and a bijection that associates members of the partition with terminal nodes of the tree. For each member of any sequence of nested subtrees of $T$ with common root, it is required that there exist a bijection that associates that subtree with a corresponding subpartition, where the nesting of partitions and of subtrees correspond in an obvious way. A real-valued function $h$ on $\Omega$ is a *binary tree-structured function* if there is a finite, rooted, binary tree-structured partition of $\Omega$ and $h$ is constant on each member of the partition.

Throughout, $d$-dimensional Euclidean space is denoted by $\mathcal{R}^d$, an important subset being the open unit cube $\mathcal{U}^d$. Much mathematics concerns the case $d = 2$. Our principal focus is on rooted, binary tree-structured partitions of $\mathcal{R}^d$ (or $\mathcal{U}^d$) into *boxes*.

A *box* is a set $\{\mathbf{x}\} = B \subset \mathcal{R}^d$ that is the solution set of a system of inequalities defined by inner products $\mathbf{b}^{(k)} \cdot \mathbf{x} \le c$ or $\mathbf{b}^{(k)} \cdot \mathbf{x} > c$, $k = 1, 2, \ldots, K < \infty$, where $\mathcal{R}^d \ni \mathbf{b}^{(k)} \ne 0$ and $c$ is real. If for each linear inequality that defines $B$ exactly one coordinate is not 0, then $B$ is a *basic box* or, alternatively, an *interval*.

Our focus is on rooted, binary tree-structured partitions of $\mathcal{R}^d$ (alternatively, of $\mathcal{U}^d$) into a finite number of basic boxes. There is an obvious bijection that associates terminal leaves of the tree and basic boxes of the partition without nonempty subsets that are themselves basic boxes of the partition. $Q$ is a generic symbol for a finite partition of $\mathcal{R}^d$ (or $\mathcal{U}^d$), all of whose component subsets are basic boxes. For $\mathbf{x} \in \mathcal{R}^d$, denote by $B(\mathbf{x})$ the unique, smallest basic box in $Q$ that contains $\mathbf{x}$. For a sequence of such partitions $Q^{(N)}$, $B^{(N)}(\mathbf{x})$ has an obvious meaning. Write $B^{(N)}$ for $\{B^{(N)}(\mathbf{x})\}$.

The *diameter* of $B^N(\mathbf{x})$ is defined as

$$\delta_N(\mathbf{x}) = \sup\{\|\mathbf{z} - \mathbf{y}\| : y, z \in B^{(N)}(\mathbf{x})\},$$

while the *norm* of $Q^{(N)}$, $\|Q^{(N)}\|$, is defined as

$$\|Q^{(N)}\| = \max_x \delta_N(x).$$

The assumptions entail that writing "max" makes sense because $\{\delta_N(\mathbf{x})\}$ is finite.

Suppose that $(\mathbf{X}, Y), (\mathbf{X}_1, Y_1), \ldots, (\mathbf{X}_N, Y_N)$ are independent, identically distributed (i.i.d.) vectors,

$$(1.1) \qquad \mathbf{X} \in \mathcal{R}^d, \ Y \in \mathcal{R}^1, \ E(|Y|) < \infty.$$

Write

$$(1.2) \qquad h(\mathbf{x}) = E(Y | \mathbf{X} = \mathbf{x})$$

for the regression of $Y$ on $\mathbf{X}$. The *test case* $\mathbf{X}$ and *learning sample* $\{(\mathbf{X}_i, Y_i) : i = 1, \ldots, N\}$ are given; $h$ is to be estimated. Write

$$(1.3) \qquad \hat{h}_N = \hat{h}_N(\mathbf{X}) = \hat{h}_N(\mathbf{X}, (\mathbf{X}_1, Y_1), \ldots, (\mathbf{X}_N, Y_N)).$$

Throughout, $\hat{h}_N$ is a simple average of those $Y_i$'s, $1 \le i \le N$, for which $\mathbf{X}_i$ lies in the same box of a partition of $\mathcal{R}^d$ as $\mathbf{X}$, provided the box has positive empirical probability. Equalities (1.3) are made precise in what follows by (1.7) and (1.8).

For a measurable subset $S \subset \mathcal{R}^d$, *define* $\mu(S)$ and $F(S)$ as

$$(1.4) \qquad \mu(S) = E(Y I_S(\mathbf{X})), \qquad F(S) = E(I_S(\mathbf{X})) = P(\mathbf{X} \in S).$$

$I_S(\mathbf{X}) = 1$ if $\mathbf{X} \in S$ and is 0 otherwise. The reader can check that a version of $h(\mathbf{x})$ is $\frac{d\mu}{dF}(x)$. For $P(\mathbf{X} \in B^{(N)}(\mathbf{x})) > 0$, *define*

$$(1.5) \qquad h_N(\mathbf{x}) = \frac{\mu(B^{(N)}(\mathbf{x}))}{P(\mathbf{X} \in B^{(N)}(\mathbf{x}))}.$$

Of course, we do not observe $h_N$ in applications. At each stage $N$ of sampling, we are given $Q^{(N)}$, a finite, rooted, binary tree-structured partition of $\mathcal{R}^d$ into basic boxes that depends measurably on $\{(\mathbf{X}_n, Y_n) : n = 1, \ldots, N\}$. Define $\mathcal{F}_N$ by

$$(1.6) \qquad \mathcal{F}_N = \sigma\{(\mathbf{X}_n, Y_n) : n = 1, \ldots, N : I_B(\mathbf{X}) : B \in Q^{(N)}\},$$

where $\sigma\{\cdot\}$ is the $\sigma$-field generated by the random quantities inside $\{\cdot\}$. I quote a lemma that appears as Lemma 3.12 in [3].

LEMMA 1.1.  $E(Y|\mathcal{F}_N) = h_N(\mathbf{X})$.

DEFINITIONS AND NOTATION.  For $\mathbf{x} \in \mathcal{R}^d$ and $B^{(N)}(\mathbf{x}) = B$, write

$$(1.7) \qquad \hat{h}_N(\mathbf{x}) = \frac{\hat{\mu}_N(B)}{\hat{F}_N(B)} I_{\{\hat{F}_N(B) > 0\}},$$

where

$$(1.8) \qquad \hat{\mu}_N(B) = \frac{1}{N} \sum_{n=1}^{N} Y_n I_B(\mathbf{X}_n)$$

and $\hat{F}_N(B) = \frac{1}{N} \sum_{n=1}^{N} I_B(\mathbf{X}_n)$.

Note from Lemma 1.1 and (1.7) and (1.8) that $\hat{h}_N(\mathbf{x})$ bears the same relationship to $B^{(N)}$ and $\hat{F}_N$ that $h_N(\mathbf{x})$ does to $B^{(N)}$ and $F$.

One can show that $\hat{h}_N(\mathbf{X}) \to h(\mathbf{X})$ in various senses as $N$ grows without bound; see, for example, [8], Theorem 12.7 of [1] and [4, 5, 10, 11, 14]. A particularly strong notion of convergence, but one that matters for applications, is unconditional almost sure convergence, where "unconditional" is meant with respect to the learning sample and test case. The question arises as to whether $\hat{h}_N$ is *consistent* in this very strong sense. A major point of this paper is that this strong notion of consistency does not generally hold.

THEOREM 1.2.  *There exists a sequence $Q^{(N)}$ of finite, rooted, binary tree-structured partitions of $\mathcal{U}^2$ for which $\|Q^{(N)}\| \to 0$, as well as a set $\{(\mathbf{X}_n, Y_n) : n = 1, 2, \ldots\}$ and an $\mathbf{X}$ that satisfy (1.1) where $\mathbf{X}$ is uniformly distributed on $\mathcal{U}^2$, $E(|Y|)$ is finite, $P(\hat{h}_N(\mathbf{X}) - h_N(\mathbf{X}) \to 0) = 1$, and yet $P(h_N(\mathbf{X}) - h(\mathbf{X}) \to 0) = 0$. Thus, the analogue of "variance" tends to 0 almost surely and the diameters of basic boxes of partitions tend to 0 surely.*

*However, the analogue of bias almost surely does not tend to* 0 *and* $P(\hat{h}_N(\mathbf{X}) \to h(\mathbf{X})) = 0$.

The perverse behavior of $h_N(\mathbf{X})$ in the example of interest here is summarized in the next theorem. Without loss of generality, we assume that $Y \geq 0$, so $h(\mathbf{x}) \geq 0$.

THEOREM 1.3. *With* $(\mathbf{X}, Y)$ *and* $\{Q^{(N)}\}$ *as in Theorem* 1.2 *and* $Y \geq 0$, $\{h : E(h(\mathbf{X})) < \infty \text{ and } \overline{\lim} h_N(\mathbf{x}) < \infty \text{ for some } \mathbf{x} \in \mathcal{U}^2\}$ *is of the first category in* $\mathcal{L}^1(\mathcal{U}^2)$.

The main examples are also relevant to understanding the (two-class) classification problem. Thus, let $Y = 1$ or 2 with probability $1/2$ each. Scale the nonnegative $h$ of Theorem 1.2 to have integral 1 and thus to be a probability density on $\mathcal{U}^2$. Suppose that given $Y = 1$, $\mathbf{X}$ has density $h$ and given $Y = 2$, $\mathbf{X}$ has the uniform distribution on $\mathcal{U}^2$. Given the training sample described in (1.1), a (measurable) empirical *classification rule* $d_N(\mathbf{x}) = d_N(\{(\mathbf{X}_i, Y_i) : i = 1, \ldots, N\})(\mathbf{x})$ is given. Then $d_N(\mathbf{x})$ takes values 1 or 2 and is a guess of the unknown $Y$ when $\mathbf{X} = \mathbf{x}$. We lose one dollar for an incorrect guess, otherwise we lose nothing. It is not difficult to show that the rule "$d_B(\mathbf{x}) = 1$ if $h(\mathbf{x}) > 1$, otherwise $d_B(\mathbf{x}) = 2$" is a Bayes rule. Its expected loss is the "Bayes risk" of the Bayes rule $d_B$. In practice, we would not know $h$, so we could not compute a Bayes rule.

A sequence $d_N$ of classification rules is said to be *Bayes-risk consistent* if the sequence of expected losses converges to the Bayes-risk of a Bayes rule as $N$ increases without bound.

COROLLARY 1.4. *For any* $\varepsilon > 0$ *and the stated problem of two-class classification, with the learning sample as in* (1.1) *and that which precedes it and with* $\{Q^{(N)}\}$ *as in Theorem* 1.2, *there is a sequence of rooted, binary tree-structured classification rules* $d_N$ *with the following property: the rules are Bayes-risk consistent, but* $P(d_N(\mathbf{X}) \to d_B(\mathbf{X})) < \varepsilon$. *As before,* $d_B$ *is a Bayes rule for the problem.*

**2. Martingales and the differentiation of integrals.** The goal of this section is to lend perspective to the main examples of this paper.

Given an $\mathcal{L}^1$ function $f$ on a probability space and a monotonic sequence of sub-$\sigma$-fields of a base $\sigma$-field, the martingale convergence theorem ensures that the sequence of successive conditional expectations converges almost surely to the conditional expectation given the "limit" $\sigma$-field. We are interested in the case where the sequence is monotonically increasing. If $f$ is measurable with respect to the $\sigma$-field ($\mathcal{T}$, say) that is generated by the sequence, then the limit random variable is $f$ itself, at least up to a $\mathcal{T}$ set of probability 0.

With our notions of basic box and $\hat{h}_N(\mathbf{X})$, provided each hyperplane determined by the boundary of each $B_m \in Q^{(N)}$ contains at least one $\mathbf{X}_n, 1 \leq n \leq N$, then $h_N$ is equivariant to strictly monotonic transformations of the coordinate axes. Thus, if $T : \mathcal{R}^d \to \mathcal{R}^d$ is of the form $T(x_1, \ldots, x_d) = (h_1(x_1), \ldots, h_d(x_d))$ with $h_i$ strictly monotone, then $T$ maps basic boxes to basic boxes—in an abuse of notation, $\hat{h}_N(T(\mathbf{x})) \equiv \hat{h}_N(\mathbf{x})$. Note that nothing is lost if we take the range of $\mathbf{X}$ to be $\mathcal{U}^d$. Clearly, mappings such as $T$ do not preserve ratios of sides of boxes.

Our application allows restrictions on neighborhood systems that are different from bounds on the ratios of sides of boxes. For one, rectangular neighborhoods are always members of finite partitions of $\mathcal{U}^2$, each of which is a rooted, binary tree-structured partition. The members of a partition are the atoms of a finite $\sigma$-field of subsets of $\mathcal{U}^2$. Because the finite $\sigma$-fields are shown not to differentiate every $\mathcal{L}^1$ function, or even "most" such functions, they cannot be nested. Even so, Gordon and Olshen ([8], Section 6) asked if the restrictions to the particularly simple probability space and shapes of the atoms would allow the relaxation of assumptions on nesting and thereby an extension of the martingale theorem, not to mention extensions of theorems regarding the almost sure consistency of binary tree-structured algorithms for regression.

DEFINITION.   For $B \subset \mathcal{R}^n$, let

$$\delta(B) = \sup\{\|\mathbf{z} - \mathbf{y}\| : \mathbf{y}, \mathbf{z} \in B\}.$$

For each $\mathbf{x} \in \mathcal{R}^n$, suppose that $\mathcal{B}(\mathbf{x})$ is a collection of bounded Borel sets with positive Lebesgue measure, that $B \in \mathcal{B}(\mathbf{x})$ implies $\mathbf{x} \in B$ and that for each $\mathbf{x}$, there exists a sequence $R_k = R_k(\mathbf{x}) \subset \mathcal{B}(\mathbf{x})$ for which $\delta(R_k) \to 0$ as $k \to \infty$. Then

$$\mathcal{B} = \bigcup_{\mathbf{x} \in \mathcal{R}^n} \mathcal{B}(\mathbf{x}) \qquad \text{is a } differentiation \ basis.$$

Write $f \in \mathcal{L}^1(\mathcal{R}^n)$ if a version of $f$ is Borel and the Lebesgue integral $\int_{\mathcal{R}^n} |f(\mathbf{x})| \lambda(d\mathbf{x})$ is finite. If $f \in \mathcal{L}^1(\mathcal{R}^n)$ and $\mathcal{B}$ is a differentiation basis, then $\mathcal{B}$ *differentiates* $\mathcal{L}^1(\mathcal{R}^n)$ if for Lebesgue almost every $\mathbf{x}$, $\mathbf{x} \in B_k$, $k = 1, 2, \ldots$, $B_k \in \mathcal{B}$ and $\delta(B_k) \to 0$ implies

$$\lim_{k \to \infty} (\lambda(B_k))^{-1} \int_{B_k} f(\mathbf{u}) \lambda(d\mathbf{u}) = f(\mathbf{x}).$$

Write $\mathcal{B}_1(\mathbf{x})$ for the collection of open, bounded cubes containing $\mathbf{x} \in \mathcal{R}^n$ and $\mathcal{B}_1^*(\mathbf{x})$ for the collection of open, bounded cubes centered at $\mathbf{x}$. Write $\mathcal{B}_1 = \bigcup_{\mathbf{x}} \mathcal{B}_1(\mathbf{x})$ and $\mathcal{B}_1^* = \bigcup_{\mathbf{x}} \mathcal{B}_1^*(\mathbf{x})$. Both $\mathcal{B}_1$ and $\mathcal{B}_1^*$ are differentiation bases. The Lebesgue differentiation theorem says that $\mathcal{B}_1^*$ differentiates $\mathcal{L}^1(\mathcal{R}^n)$. Also, $\mathcal{B}_1$ differentiates $\mathcal{L}^1$. These conclusions remain true if the definitions of $\mathcal{B}_1$ and $\mathcal{B}_1^*$ are relaxed to allow their members to be basic boxes instead of cubic intervals, but with a finite bound on the ratio of dimensions of any two sides of the boxes.

Let $\mathcal{B}_2(\mathbf{x})$ be the set of otherwise unrestricted basic boxes that contain $\mathbf{x} \in \mathcal{R}^n$, and define $\mathcal{B}_2^*(\mathbf{x})$ by analogy. Let $\mathcal{B}_2 = \bigcup_{\mathbf{x}} \mathcal{B}_2(\mathbf{x})$ and $\mathcal{B}_2^* = \bigcup_{\mathbf{x}} \mathcal{B}_2^*(\mathbf{x})$. Neither $\mathcal{B}_2$ nor even $\mathcal{B}_2^*$ differentiates $\mathcal{L}^1(\mathcal{R}^n)$ or, for that matter, $\mathcal{L}^1(\mathcal{U}^2)$ (pages 95 and 96 of [8]). If in addition to $f \in \mathcal{L}^1(\mathcal{R}^n)$, we also have

$$\int |f(\mathbf{x})|(1 + \max(0, \log |f(\mathbf{x})|))^{n-1}\lambda(d\mathbf{x}) < \infty,$$

then $\mathcal{B}_2$ differentiates $f$.

The author believes the proof of Lemma 2.1 given below to be new. This lemma is at the heart of the counterexamples.

LEMMA 2.1. *For $N = 3, 4, 5, \ldots$, there is a nonnegative Borel $f_N$ on $\mathcal{U}^2$ for which:*

(i) *$\int_{\mathcal{U}^2} f_N(\mathbf{u})\lambda(d\mathbf{u}) \leq N^{-1}$;*

(ii) *for each $\mathbf{x} \in \mathcal{U}^2$, there exists a basic box $B_{2,N}(\mathbf{x})$ with $\delta(B_{2,N}(\mathbf{x})) < N^{-1}$;*

(iii) *$\int_{B_{2,N}(\mathbf{x})} f_N(\mathbf{u})\lambda(d\mathbf{u}) > N\lambda(B_{2,N}(\mathbf{x}))$.*

Obviously, $\{B_{2,N}(\mathbf{x})\}$ can here be taken to be open intervals.

PROOF OF LEMMA 2.1. Let $\varepsilon > 0$ be given. Further, set $\beta_n = 1/n(\ln n)^2$, $\eta_n = (\ln n)^{-1}$ and $\gamma_n = (\ln n)^{1/2}$. Write $S_n$ for $\{0 < x \leq \beta_n/\eta_n, 0 < y \leq \eta_n\} \cup \{0 < x \leq \eta_n, 0 < y \leq \beta_n/\eta_n\}$. The set $S_n$ is the union of two oblong rectangles, one contiguous to the $x$-axis and one contiguous to the $y$-axis. For $(x, y) = \mathbf{u} \in \mathcal{U}^2$, $g_n = g_n(\mathbf{u})$ is defined to be $\gamma_n I_{S_n}(\mathbf{u})$. Obviously, $g_n \geq 0$ and $\gamma_n \beta_n \leq \|g_n\|_1 \leq 2\gamma_n \beta_n$. Write $R_n = S_n \cup \{\mathbf{u} \in \mathcal{U}^2 : xy \leq \beta_n, (\beta_n/\eta_n) < x < \eta_n, (\beta_n/\eta_n) < y < \eta_n\}$. Thus, $R_n$ is the union of $S_n$ and a set that is bounded by the $x$-axis, the $y$-axis, $\{x = \eta_n\}$, $\{y = \eta_n\}$ and the hyperbola $\{xy = \beta_n\}$. Furthermore, the hyperbola has nonempty intersection with $S_n$. Note that for $\mathbf{u} \in R_n$, there exists a basic box $R_n'(\mathbf{u}) = R_n' \subset R_n$ bounded by the $x$-axis, the $y$-axis and with a vertex on the hyperbola $\{xy = \beta_n\}$ such that

$$\int_{R_n'} g_n/\lambda(R_n') \geq \frac{1}{2}\gamma_n\beta_n\beta_n = \frac{1}{2}\gamma_n \nearrow \infty.$$

$\{R'_n(\mathbf{u})\}$ can be assumed to consist of only a countable class of open subsets of $R_n$. Also,

$$\lambda(R_n) = \beta_n + \beta_n \int_{\beta_n/\eta_n}^{\eta_n} dx/x = \beta_n + 2\beta_n \ln \eta_n + (-\beta_n \ln \beta_n).$$

Therefore, $\sum \lambda(R_n) = \infty$. On the other hand, because $\sum \gamma_n \beta_n$ converges, there exists an $N = N(\varepsilon)$ sufficiently large that $\sum_N^\infty \gamma_n \beta_n < \varepsilon/2$, so $\| \sum_N^\infty g_n \|_1 < \varepsilon$.

For $n = N, N+1, \ldots$ on the square with vertices $(\eta_n, \eta_n)$, $(\eta_n, 1 - \eta_n)$, $(1 - \eta_n, \eta_n)$, $(1 - \eta_n, 1 - \eta_n)$, choose a point $P_n$ uniformly at random so that $P_N, P_{N+1}, \ldots$ are independent. Place a square with sides $\eta_n$ in the cited (larger) square so that the center is at $P_n$ and the sides are parallel to the coordinate axes. Call this random square $(\mathrm{Sq})_n$. Three subsets of $(\mathrm{Sq})_n$ require definition.

In what follows, two planar sets are *homothetic* if one is identical to the other up to a rigid motion of the plane not involving rotation. Denote by $\mathcal{S}_n$ the subset of $(\mathrm{Sq})_n$ that is homothetic to $S_n$, by $\mathcal{R}_n$ the subset of $(\mathrm{Sq})_n$ that is homothetic to $R_n$ and by $\mathcal{R}'_n$ the subset of $(\mathrm{Sq})_n$ that is a basic box and is homothetic to the basic box $R'_n$. For $\mathbf{u} \in \mathcal{U}^2$ define the (random) function $h = h(\mathbf{u})$ by $h = \sum_N^\infty \gamma_n I_{\mathcal{S}_n}$. Necessarily, $0 \le h$ and $\|h\|_1 < \varepsilon$. Because $\sum \lambda(R_n) = \infty$, manipulation of indicator functions, independence and monotone convergence guarantee that almost surely

$$\lambda\left(\bigcup \mathcal{R}_n\right) = 1.$$

Moreover,

$$\int_{\mathcal{R}'_n} h/\lambda(\mathcal{R}'_n) \ge \frac{1}{2} \gamma_n \nearrow \infty.$$

It now follows from Fubini's theorem that there exists a real-valued function $g$ on $\mathcal{U}^2$ for which (i) $0 \le g \in \mathcal{L}^1(\mathcal{U}^2)$, (ii) $\|g\|_1 < \varepsilon$ and (iii) for almost all $\mathbf{u}_0 \in \mathcal{U}^2$, there exists a basic box $\mathcal{R}''_n = \mathcal{R}''_n(\mathbf{u}_0) \subset \mathcal{U}^2$ with $\mathbf{u}_0 \in \mathcal{R}''_n$ and

$$(2.1) \qquad\qquad \int_{\mathcal{R}''_n} g/\lambda(\mathcal{R}''_n) \ge \varepsilon^{-1},$$

at least for $n$ satisfying $\frac{1}{2}\gamma_n > \varepsilon^{-1}$.

Finally, one can argue that "almost all $\mathbf{u}_0$" can, in fact, be "all." With the $\mathcal{R}''_n$ open, the set of $\mathbf{u}_0$ of full measure for which (2.1) holds is seen to be open. Its complement is thus a closed set of Lebesgue measure 0. Call it $N$. One sees that there exists $G \in \mathcal{L}^1(\mathcal{U}^2)$ that is continuous on $\mathcal{U}^2 \backslash N$ and that tends to $\infty$ as its argument tends to $N$. Without loss of generality, $G \ge 0$. Now let $f = g + G$ in order to establish Lemma 2.1. $\square$

Definition. For $Q$ a finite partition of $\mathcal{U}^2$ into basic boxes $B$ with $\lambda(B) > 0$ for all $B \in Q$, $f \in \mathcal{L}^1(\mathcal{U}^2)$ and $\mathbf{u} \in \mathcal{U}^2$, define

$$(2.2) \qquad E(f|Q)(\mathbf{u}) = \sum_{B \in Q} I_B(\mathbf{u}) \left( \int_B f(\mathbf{x}) \lambda(d\mathbf{x}) / \lambda(B) \right).$$

**3. Examples.** Material in this section expands upon that of Lemma 2.1 and is at the heart of the proofs of Theorems 1.2 and 1.3 and Corollary 1.4. Compare the results here with those of Busemann and Feller [2] and also with those of Saks [13].

*Define* $\mathcal{U}_N^2$ to be $\{\mathbf{u} \in \mathcal{U}^2, \mathbf{u} = (s,t) : N^{-1} \le s, t \le 1 - N^{-1}\}$. Therefore, the (countable class of) open basic boxes $\{\mathcal{R}_n''(\mathbf{u}) : \mathbf{u} \in \mathcal{U}_N^2\}$, whose existence is ensured by Lemma 2.1, is an open cover of $\mathcal{U}_N^2$. Because $\mathcal{U}_N^2$ is compact in the usual topology, the Heine–Borel theorem guarantees the existence of a finite subcover of open basic boxes, which we denote by $\{B_{2,N}(\mathbf{u}_j) : j = 1, \ldots, K_N\}$. Now, fix $N$, $j$, $1 \le j \le N$, and $B_{2,N}(\mathbf{u}_j)$. There is clearly a rooted, finite, binary tree-structured partition $Q^{(N,j)}$ of $\mathcal{U}^2$ for which $B_{2,N}(\mathbf{u}_j) \in Q^{(N,j)}$, $\|Q^{(N,j)}\| < N^{-1}$ and $E(f_N|Q^{(N,j)})(\mathbf{u}) > N$ for $\mathbf{u} \in B_{2,N}(\mathbf{u}_j)$. We therefore have the following lemma:

Lemma 3.1. *For $j = 1, \ldots, K_N$, $\|Q^{(N,j)}\| < N^{-1}$ and for all $\mathbf{x} \in \mathcal{U}_N^2$,*

$$\max\{E(f_N|Q^{(N,j)})(\mathbf{u}) : j = 1, \ldots, K_N\} > N.$$

Now define $f = f(\mathbf{u})$ on $\mathcal{U}^2$ by

$$(3.1) \qquad f = \sum_{N=3}^{\infty} f_{N^2}.$$

Clearly, $f \in \mathcal{L}^1(\mathcal{U}^2)$. The finite, rooted binary tree-structured partitions of $\mathcal{U}^2$ that concern us are

$$(3.2) \qquad \ldots, Q^{(N^2,1)}, Q^{(N^2,2)}, Q^{(N^2,K_{N^2})}, \ldots, Q^{((N+1)^2,1)}, \ldots,$$
$$Q^{((N+1)^2, K_{(N+1)^2})}, \ldots.$$

From Lemma 3.1, it follows that $\|Q^{(N^2,j)}\| < N^{-2}$. Since $f \ge f_{N^2}$ ($N = 3, 4, \ldots$), it follows also from that lemma that

$$\max\{E(f|Q^{(N^2,1)})(\mathbf{u}), \ldots, E(f|Q^{(N^2,K_{N^2})})(\mathbf{u})\} > N^2$$

if $\mathbf{u} \in \mathcal{U}_{N^2}^2$. Therefore, if we relabel the partitions (3.2) as

$$(3.3) \qquad \ldots, Q^{(n)}, Q^{(n+1)}, \ldots,$$

then Theorem 3.2 follows.

THEOREM 3.2.  *With $Q^{(n)}$ as in (3.3), as $n$ grows without bound, $\|Q^{(n)}\| \to$ $0$, but $\overline{\lim} E(f|Q^{(n)})(\mathbf{u}) = +\infty$ on $\mathcal{U}^2$.*

We continue now with a corollary that is due to R.M. Dudley and to Gordon and Olshen; see ([3], pp. 161–162).

COROLLARY 3.3.  *There exist a probability $\nu$ on $\mathcal{U}$ and an open set $\mathcal{O} \subset U$ with the property that*

$$\nu\{\overline{\lim} E(I_\mathcal{O}|Q^{(n)}) > 0\} > \nu(\mathcal{O}).$$

*Here, $Q^{(n)}$ is as in (3.3).*

Next, we argue that the analogue of Theorem 1.3 is true in the present context. First, note that $\{Q^{(n)}\}$ of (3.3) can be defined inductively so that for each $\mathbf{u} \in \mathcal{U}^2$, $|\bigcup_{n=3}^\infty \partial(B^{(n)}(\mathbf{u}))| \leq 4$. We will assume this to be the case. Here, for $S \subset \mathcal{U}^2$, $\partial(S)$ denotes its boundary and $|S|$ its cardinality.

THEOREM 3.4.  *Let $F = \{g \in \mathcal{L}^1(\mathcal{U}^2) : \overline{\lim} E(g|Q^{(n)}(\mathbf{u})) < \infty$, some $\mathbf{u} \in \mathcal{U}^2\}$. Then $F$ is of the first category in $\mathcal{L}^1(\mathcal{U}^2)$.*

PROOF.  For $k = 1, 2, \ldots$ and $M = 3, 4, \ldots$, let $F_{k,M} = \{f \in \mathcal{L}^1(\mathcal{U}^2):$ for some $\mathbf{u} \in \mathcal{U}_M^2, \|B^{(n)}(\mathbf{u})\| \leq k^{-1}$ implies $E(f|Q^{(n)})(\mathbf{u}) \leq k\}$. Necessarily, $F = \bigcup_{k,M} F_{k,M}$. Therefore, it is enough to show that each $F_{k,M}$ is of the first category in $\mathcal{L}^1(\mathcal{U}^2)$. To that end, fix $k, M$ and assume that $\{g_{j,k,M}\}$ (each in $F_{k,M}$) satisfies $E(|g_{j,k,M} - g|) \to 0$ as $j \to \infty$ for some $g \in \mathcal{L}^1(\mathcal{U}^2)$. By the definition of $F_{k,M}$, for each $j = 1, 2, \ldots$, there exists $\mathbf{u}_{j,k,M} \in \mathcal{U}_M^2$ for which $\|B^{(n)}(\mathbf{u}_{j,k,M})\| \leq k^{-1}$ implies that $E(g_{j,k,M}|Q^{(n)})(\mathbf{u}_{j,k,M}) \leq k$. $\mathcal{U}_M^2$ is compact and without loss of generality, we can assume that $\mathbf{u}_{j,k,M} \to \mathbf{u}_{k,M}$ as $j \to \infty$. For $n = n(\mathbf{u}_{k,M})$ sufficiently large, $\mathbf{u}_{k,M} \notin \bigcup_{l=n}^\infty \partial(B^{(l)}(\mathbf{u}_{k,M}))$. For $l$ sufficiently large, $\|B^{(l)}(\mathbf{u}_{k,M})\| \leq k^{-1}$. Eventually, for each fixed $l$, $\mathbf{u}_{j,k,M} \in B^{(l)}(\mathbf{u}_{k,M})$. Therefore, $B^{(l)}(\mathbf{u}_{k,M}) = B^{(l)}(\mathbf{u}_{j,k,M})$ for $j$ sufficiently large. For such $j$, $E(g_{j,k,M}|Q^{(l)})(\mathbf{u}_{k,M}) \leq k$. Now, $E(g|Q^{(l)})(\mathbf{u}_{k,M}) \leq |E(g - g_{j,k,M}|Q^{(l)})(\mathbf{u}_{k,M})| + k \leq E(|g - g_{j,k,M}||Q^{(l)})(\mathbf{u}_{k,M}) + k$. By making $l$, then $j$, sufficiently large, we can conclude that $g \in F_{k,M}$, that is, that $F_{k,M}$ is closed.

We now show that $F_{k,M}$ contains no open ball in $\mathcal{L}^1(\mathcal{U}^2)$. Let $h \in F_{k,M}$, $f$ be as in the example and $0 < \alpha < 1$. Because every function in $\mathcal{L}^1$ is the limit of bounded, continuous functions and $g$ is taken to be in an open $\mathcal{L}^1$ ball in $F_{k,M}$, we can (and do) take $h$ to be bounded and continuous. Define $g_\alpha = (1 - \alpha)h + \alpha f$. Then $E(|h - g_\alpha|) \to 0$ as $\alpha \to 0$. But for each fixed $\alpha$, $\overline{\lim} E(g_\alpha|Q^{(n)}) = \infty$ everywhere on $\mathcal{U}^2$. This completes the proof of Theorem 3.4.  □

**4. Applications to rooted, binary tree-structured regression and classification.** We return, now, in this last section, to the application of the results of Section 3 to arguments for Theorems 1.2 and 1.3 and Corollary 1.4. The point of Theorem 1.2 is that there exists a sequence $Q^{(N)}$ of finite, rooted, binary tree-structured partitions of the unit cube $\mathcal{U}^2$ in $\mathcal{R}^2$ for which $\|Q^{(N)}\| \to 0$, as well as a set $\{(\mathbf{X}_n, Y_n) : n = 1, 2, \ldots\}$ and an $\mathbf{X}$ that satisfy the assumptions given previously, for which $\mathbf{X}$ is uniformly distributed on $\mathcal{U}^2$ and $E(|Y|)$ is finite, but where $P(\hat{h}_N(\mathbf{X}) \to h(\mathbf{X})) = 0$. Write

$$(4.1) \quad |\hat{h}_N(\mathbf{X}) - h(\mathbf{X})| \geq |h_N(\mathbf{X}) - h(\mathbf{X})| - |\hat{h}_N(\mathbf{X}) - h_N(\mathbf{X})| := II - I.$$

The original question posed by Gordon and Olshen pertained both to $|\hat{h}_N(\mathbf{X}) - h(\mathbf{X})|$ and to $II$. If the counterexample to the almost everywhere convergence to 0 of $II$ for an $h$ in $\mathcal{L}^1(\mathcal{U}^2)$ implies the existence of an analogous counterexample to the almost sure convergence of $|\hat{h}_N(\mathbf{X}) - h(\mathbf{X})|$ to 0 for an $h$ with $E(|h(\mathbf{X})|) < \infty$, then Theorem 1.2 is proved. In fact, $I$ can converge to 0 almost surely while $II$ does not.

The asymptotic behavior of $I$ depends on the large deviation behavior of

$$\sup_{D \in \mathcal{D}} |\hat{F}_N(D) - F(D)|,$$

where $\mathcal{D}$ is a Vapnik–Chervonenkis class, as is the set of basic boxes, that is, the set of interval subsets of $\mathcal{U}$. For Theorems 1.2 and 1.3 and Corollary 1.4, we do not need Vapnik–Chervonenkis-like results. It is easy to adapt Theorem 3.2 so as to preserve the lack of convergence of $II$, while in, fact, $I$ tends to 0 almost surely. Suppose that $(\mathbf{X}, Y)$ and the training sample $(\mathbf{X}_1, Y_1), \ldots, (\mathbf{X}_N, Y_N)$ are as in (1.1), with $\mathbf{X}$ having a uniform distribution on $\mathcal{U}^2$. Fix an $n$ and therefore $Q^{(n)}$ as in (3.3) and let $\mathcal{F}_n$ be as in (1.6). For $i = 1, 2, \ldots$, let $Y_i = f(\mathbf{X}_i)$, where $f$ is as in (3.1) and Theorem 3.2. Recall that $B \in Q^{(n)}$ implies that $\lambda(B) > 0$ and that each $Q^{(N)}$ has only finitely many members. Therefore, the strong law of large numbers implies that for any fixed $Q^{(n)}$ and $\varepsilon = \varepsilon_n$, $0 < \varepsilon_n < 1$, there exists an $N(n, \varepsilon_n)$ sufficiently large that $P(\bigcup_{N > N(n,\varepsilon_n)} |\hat{h}_N(\mathbf{X}) - h_N(\mathbf{X})| > \varepsilon_n) < \varepsilon_n$. If $\varepsilon_n$ is the term of a convergent series, then $N(n, \varepsilon_n)$ can grow sufficiently fast and $n$ sufficiently slowly so that $Q^{(n)}$ applies to learning samples from size $N(n, \varepsilon_n)$ to $N(n + 1, \varepsilon_{n+1})$. The Borel–Cantelli lemma implies that $I$ tends to 0 almost surely. That is, the cardinality of the learning sample, $N$, and the $Q^{(N)}$ that applies need not be related, other than for convenience. It follows that $II$ can fail to converge to 0 on a set of probability 1, while $I$ can converge to 0 with probability 1. When this happens, $|\hat{h}_N(\mathbf{X}) - h(\mathbf{X})|$ fails to converge to 0 on a set of probability 1. This completes the argument for Theorem 1.2.

A repetition of the argument in the previous paragraph, with Theorem 3.4 substituted for Theorem 3.2, completes the argument for Theorem 1.3.

Because $h = h(\mathbf{u})$ can be approximated arbitrarily closely in $\mathcal{L}^1(\mathcal{U}^2)$ by a continuous function, it is clear that in the extended example of this section, $\hat{h}_N(\mathbf{X}) - h(\mathbf{X})$ tends to 0 in $\mathcal{L}^1$ of the common probability space on which random variables $(\mathbf{X}, Y)$ and the learning sample are defined. This observation is analogous to the argument for Proposition 1 of [15].

We now turn our attention to a brief discussion of the two-class classification problem and argument for Corollary 1.4. A formulation is given after Theorem 1.3 of Section 1. We argue that the rule for two-class classification given next is Bayes-risk consistent, but not almost surely consistent.

With $Q^{(N)}$ as in the arguments for Theorems 1.2 and 1.3, let $d_N(\mathbf{x}) = 1$ if

$$\sum_{B \in Q^{(N)}} \sum_{i=1}^N I_{[\mathbf{x} \in B, \mathbf{X}_i \in B, Y=1]} > \sum_{B \in Q^{(N)}} \sum_{i=1}^N I_{[\mathbf{x} \in B, \mathbf{X}_i \in B, Y=2]},$$

otherwise, $d_N(\mathbf{x}) = 2$. It follows from the construction of $\{Q^{(N)}\}$ and Theorem 12.17 of [1] that $d_N$ is Bayes-risk consistent. From the argument for Theorem 1.2 in this section, it follows that for any $\varepsilon > 0$, $Q^{(N)}$ and $h$ can be chosen so that $P(h(\mathbf{X}) < 1) > 1 - \varepsilon$, but $P(h(\mathbf{X}) < 1; d_N \to d_B) = 0$.

**Acknowledgments.** Work on this paper began long ago in the form of conversations with Lou Gordon and also with Dick Dudley. More recently, Persi Diaconis encouraged discussion among Alexandra Bellow, Izzy Katznelson, Don Ornstein and me on many aspects of the research. An Editor, an Associate Editor and a referee gave good advice. Bonnie Chung helped with technical matters. While I am grateful for the considerable help, those who supplied it bear no responsibility for any flaws that remain. Izzy Katznelson showed me why the set of measure 0 that appears at the end of the argument for Lemma 2.1 can be taken to be empty. I learned the construction in the proof of Theorem 3.2 from Alexandra Bellow.

## REFERENCES

[1] Breiman, L., Friedman, J. H., Olshen, R. A. and Stone, C. J. (1984). *Classification and Regression Trees*. Wadsworth, Belmont, CA. Since 1993 this book has been published by Chapman and Hall, New York. MR0726392

[2] Busemann, H. and Feller, W. (1934). Zur Differentiation der Lebesgueschen Integrale. *Fund. Math.* **22** 226–256.

[3] de Guzmán, M. (1975). *Differentiation of Integrals in $R^n$*. Lecture Notes in Math. **481**. Springer, Berlin. MR0457661

[4] Devroye, L., Györfi, L. and Lugosi, G. (1996). *A Probabilistic Theory of Pattern Recognition*. Springer, New York. MR1383093

[5] Devroye, L. and Krzyżak, A. (2002). New multivariate product density estimators. *J. Multivariate Anal.* **82** 88–110. MR1918616

[6] Donoho, D. L. (1997). CART and best-ortho-basis: A connection. *Ann. Statist.* **25** 1870–1911. MR1474073

[7] GERSHO, A. and GRAY, R. M. (1992). *Vector Quantization and Signal Compression.* Kluwer, Dordrecht.

[8] GORDON, L. and OLSHEN, R. A. (1984). Almost surely consistent nonparametric regression from recursive partitioning schemes. *J. Multivariate Anal.* **15** 147–163. MR0763592

[9] HASTIE, T., TIBSHIRANI, R. and FRIEDMAN, J. (2001). *The Elements of Statistical Learning: Data Mining, Inference and Prediction.* Springer, New York. MR1851606

[10] LUGOSI, G. and NOBEL, A. (1996). Consistency of data-driven histogram methods for density estimation and classification. *Ann. Statist.* **24** 687–706. MR1394983

[11] NOBEL, A. (1996). Histogram regression estimation using data-dependent partitions. *Ann. Statist.* **24** 1084–1105. MR1401839

[12] RIPLEY, B. D. (1996). *Pattern Recognition and Neural Networks.* Cambridge Univ. Press. MR1438788

[13] SAKS, S. (1934). Remarks on the differentiability of the Lebesgue indefinite integral. *Fund. Math.* **22** 257–261.

[14] STONE, C. J. (1977). Consistent nonparametric regression (with discussion). *Ann. Statist.* **5** 595–645. MR0443204

[15] ZHANG, H. and SINGER, B. (1999). *Recursive Partitioning in the Health Sciences.* Springer, New York. MR1683316

SCHOOL OF MEDICINE
HRP REDWOOD BUILDING
STANFORD UNIVERSITY
STANFORD, CALIFORNIA 94305-5405
USA
E-MAIL: olshen@stanford.edu